\RequirePackage{fix-cm}
\documentclass[smallextended]{svjour3}       % onecolumn (second format)
\smartqed  % flush right qed marks, e.g. at end of proof
\usepackage{graphicx}
\usepackage{amssymb}
\RequirePackage{color}

\def\RR{\hbox{\sf I\kern-.14em\hbox{R}}}

\begin{document}

\title{Maximum incremental (minimum decremental) ratio and one-sided RTS in FDH production technologies%\thanks{Grants or other notes
}
%\subtitle{}

\titlerunning{Maximum incremental (minimum decremental) ratios and RTS in FDH}        % if too long for running head

\author{Amin Mostafaee         \and
        Majid Soleimani-damaneh %etc.
}


\institute{Amin Mostafaee \at
              Department of Mathematics, College of Science, Tehran North Branch, Islamic Azad University, Tehran, Iran. \\
                        %  \\
           \and
           Majid Soleimani-damaneh (Corresponding author)\at
              School of Mathematics, Statistics and
Computer Science, College of Science, University of Tehran,
Enghelab Avenue, Tehran, Iran.\\
Tel.: +98-21-61112613\\
              Fax: +98-21-66412178\\
              \email{soleimani@khayam.ut.ac.ir}
}

\date{Received: date / Accepted: date}

\maketitle

\begin{abstract}
In very recent decades, Scale Elasticity (SE), as a quantitative characterization of Returns to Scale (RTS), has been an attractive research issue in convex DEA models. However, we show that the existing Scale Elasticity (SE) measure does not work properly under nonconvex FDH technologies. Due to this, we define an SE counterpart in FDH models. To this end, two new quantities, called maximum incremental and minimum
decremental ratios are introduced; and a polynomial-time procedure is developed to calculate them. The second part of the paper is devoted to introducing and investigating
left- and right-hand RTS notions in FDH models. Some necessary and
sufficient conditions are established, leading to a polynomial-time test for identification of the one-sided RTS status of DMUs. Finally, the relationships between one-sided RTS, Global RTS (GRS), and the newly-defined ratios are
established.

\end{abstract}

%============================
\section{Introduction}
 An important class of Data envelopment analysis (DEA) models is that of Free Disposal Hull (FDH) models. These models, first addressed by \textcolor{blue}{Deprins, Simar and Tulkens (1984)}, evaluate the Decision Making Units (DMUs) considering the closest inner approximation of the true strongly disposable (but possibly nonconvex) technology. FDH models have been studied by various researchers, including \textcolor{blue}{Tulkens (1993); Kerstens and Vanden Eeckaut (1999); Cherchye, Kuosmanen and Post (2001); Podinovski (2004c); Leleu (2006); Briec and Kerstens (2006); Soleimani-damaneh and Reshadi (2007); and Soleimani-damaneh and Mostafaee (2009), Kerstens and Van de Woestyne (2014)}.

Scale elasticity (SE), as a measure of response of outputs to the change(s) of inputs, is an important criterion for analyzing the performance and productivity of DMUs \textcolor{blue}{(Cooper, Seiford, and Tone, 2007).} SE can be considered as a quantitative
measure of the strength of the Returns to Scale (RTS) of DMUs. This notion has been studied by various scholars from both quantitative and qualitative standpoints; see e.g.  \textcolor{blue}{Starrett (1977); Panzar and Willig (1977); F$\ddot{a}$re, Grosskopf, and Lovell (1985); F${\o}$rsund (1996); Fukuyama (2000, 2001); Hadjicostas and Soteriou (2006); Soleimani-damaneh, Jahanshahloo, Mehrabian, and Hasannasab (2009);} and \textcolor{blue}{Podinovski and F${\o}$rsund (2010)}, among others.

Under convex technologies, SE at a given efficient DMU is defined with respect to the derivative of a Response Function (RF) which gives maximal proportion of outputs for a given proportion of inputs under feasibility; see \textcolor{blue}{Podinovski and F${\o}$rsund (2010)} and \textcolor{blue}{Podinovski, F${\o}$rsund, and Krivonozhko (2009)}. In contrast to convex technologies, in FDH ones the derivative of RF at a given efficient DMU is either zero or infinity, and so it does not provide useful information about the SE measure at the corresponding point. In this paper, we are going to extend and investigate a counterpart for SE notion in nonconvex FDH technologies. %%%We are faced with two main difficulties in this way, the nonconvexity of the FDH technologies and the non-continuity of the RF. %Convexity in DEA models helps to calculate SE by standard DEA models (see \textcolor{blue}{Fukuyama, 2000; Hadjicostas and Soteriou, 2006; and Soleimani-damaneh et al., 2009}), and continuity helps to interpret SE as a derivation of the production function (see \textcolor{blue}{Podinovski et al. 2009}). But we loss these two important properties in nonconvex FDH models.
%%%Due to this fact, we develop the concepts of
%right- and left-hand SE in such a way that these can be interpreted
%as the lower and upper estimates of elasticity measure,
%respectively.
To do this, we define two ratios, called maximum incremental and minimum decremental ratios; and we sketch a polynomial-time procedure to obtain them.

One of the concepts which is very close to SE is RTS notion. Although there are many papers for estimating RTS in FDH models (see \textcolor{blue}{Kerstens and Vanden Eeckaut, 1999; Podinovski, 2004c;
Soleimani-damaneh, Jahanshahloo, and Reshadi, 2006; Soleimani-damaneh and
Reshadi, 2007; and Soleimani-damaneh and Mostafaee, 2009}) none of the above-mentioned papers focused on one-sided RTS. In this study, we provide a PPS-based definition of right- and left-RTS under FDH technologies. A polynomial-time procedure to identify the one-sided RTS of DMUs is presented. Finally, the relationships between maximum incremental ratio, minimum decremental ratio, one-sided RTS, and Global RTS (GRS) are established.

The rest of the paper is organized as follows. Some preliminaries are given in Section 2. Failing traditional definition of SE for FDH PPSs is highlighted in Section 3. Maximum incremental and minimum decremental ratios (as SE counterparts in FDH models) are investigated in Section 4. Section 5 is devoted to defining and identifying one-sided RTS. The relationships between one-sided RTS, GRS, and above-mentioned two ratios are established in Section 6. Section 7 concludes the paper. Proofs of the main results are presented in appendix.

\section{Preliminaries}

Suppose that we have a set of $n$ DMUs consisting of DMU$_{j}$;
$j\in J=\{1,\ldots,n\}$. Each DMU$_{j}$ consumes $m$ positive inputs
$x_{1j},x_{2j},\ldots,x_{mj}$ to produce $s$ positive outputs
$y_{1j},y_{2j},\ldots,y_{sj}$. %We assume that there is not any
%duplicated DMU.
Define $x_j:=(x_{1j},x_{2j},\ldots,x_{mj})^t$ and
$y_j:=(y_{1j},y_{2j},\ldots,y_{sj})^t$ as input and output
vectors of DMU$_j$, respectively. In this paper, the superscript
$t$ stands for transpose. Also, define $X:=[x_1,x_2,\ldots,x_n]$
and $Y:=[y_1,y_2,\ldots,y_n]$ as $m\times n$ and $s\times n$
matrices of inputs and outputs, respectively.

The FDH technologies under different RTS assumptions can be represented as follows:
%--------------
\begin{eqnarray*}
T^{FDH_\Delta}=\{(x,y)~\vert~\sum_{j\in J}\lambda_jx_j \leq x,~~\sum_{j\in J} \lambda_j y_j\geq y\geq 0,~~\lambda_j=\delta \mu_j;~ \forall j\in J,\\
 \sum_{j\in J} \mu_j=1,~~\mu \in (\{0,1\})^n,~~\delta \in \Delta \}
\end{eqnarray*}
%----------------
where $\Delta$, depending on the RTS assumption of the reference
technology, is$$\Delta^{VRS}\equiv\{\delta~\vert~\delta=1\},~~~~\Delta^{CRS}\equiv\{\delta~\vert~\delta\geq 0\},$$
$$\Delta^{NIRS}\equiv\{\delta~\vert~0\leq \delta \leq 1\},~~~~\Delta^{NDRS}\equiv\{\delta~\vert~\delta\geq 1\}.$$

%-------------------------------------
Here, VRS, CRS, NIRS, and NDRS stand for variable, constant,
non-increasing, and nondecreasing RTS, respectively. Hereafter, for simplicity, we use notations FDH$_V$, FDH$_C$, FDH$_{NI}$, and FDH$_{ND}$, instead of FDH$_{VRS}$, FDH$_{CRS}$, FDH$_{NIRS}$, and FDH$_{NDRS}$, respectively.

Considering DMU$_o=(x_o,y_o);~(o \in J)$ as the unit under assessment, the
input-oriented and output-oriented FDH radial efficiency measures
of $DMU_o$ are obtained by solving the following
mixed-integer nonlinear programming problems, respectively:
\begin{equation}\label{EQNAA}
\displaystyle\begin{array}{lll}
\theta_o^\Delta=&\min& \theta\\
&s.t.&\displaystyle\sum_{j\in J} \lambda_jx_j\leq \theta x_o,\\
&&\displaystyle\sum_{j\in J} \lambda_j y_j \geq y_o,\\
&&\displaystyle\lambda_j=\delta \mu_j,~~\mu_j \in \{0,1\};~\forall j\in J,\\
&&\displaystyle\sum_{j\in J}\mu_j=1,~ \delta \in \Delta,\\
\end{array}
\end{equation}
%-----------------------------------------
\begin{equation}\label{EQNBB}
\displaystyle\begin{array}{lll}
\varphi_o^\Delta=&\max& \varphi\\
&s.t.&\displaystyle\sum_{j\in J} \lambda_jx_j\leq  x_o,\\
&&\displaystyle\sum_{j\in J} \lambda_j y_j \geq \varphi y_o,\\
&&\displaystyle\lambda_j=\delta \mu_j,~~\mu_j \in \{0,1\};~\forall j\in J,\\
&&\displaystyle\sum_{j\in J}\mu_j=1,~ \delta \in \Delta,\\
\end{array}
\end{equation}
%------------------------------
wherein $ \Delta \in \{FDH_V, FDH_C, FDH_{NI}, FDH_{ND}\}.$ %%%Due the
%%%constraints of the above models, we have
%%%\begin{equation}\label{Hba}
%%%\theta_o^{FDH_C}=\frac{1}{\varphi_o^{FDH_C}},~~~~~~~~~~~~~~~~~~~~~~~~~
%%%\end{equation}
%%%\begin{equation}
%%%\theta_o^{FDH_C}\leq \min\{\theta_o^{FDH_V}, \theta_o^{FDH_{NI}},\theta_o^{FDH_{ND}}\},
%%%\end{equation}
%%%\begin{equation}
%%%~~~\varphi_o^{FDH_C}\geq \max\{\varphi_o^{FDH_V}, \varphi_o^{FDH_{NI}},
%%%\varphi_o^{FDH_{ND}}\}.
%%%\end{equation}
%-------------------------------------
\begin{definition}
The $DMU_o=(x_o, y_o)$ is called $\Delta$-efficient if there
exists no $(x,y)\in T^{\Delta}$ such that $x\leq x_o, y\geq y_o$
and $(x,y)\neq (x_o,y_o),$ where the vector inequalities are
understood componentwise.
\end{definition}
%--------------------------
\begin{definition}\label{MPSS}
The $DMU_o=(x_o, y_o)$ is called an MPSS (i.e. a unit with Most Productive Scale Size) if $\theta_o^{FDH_C}=1$.
\end{definition}
%--------------------------

There are various papers in the DEA literature
discussing the local RTS in $\Delta-$technologies. In a local sense, DMUs are categorized to three classes: the DMUs with Increasing RTS (IRS) status, the DMUs with Decreasing RTS (DRS) status, and the DMUs with Constrant RTS (CRS) status. \textcolor{blue}{Podinovski (2004a)} pointed out that local RTS status is not a proper indicator of the direction in which MPSS in nonconvex PPSs is achieved. He defined Global RTS (GRS) by adding a fourth case called Global Sub-Constant RTS (G-SCRS) under which MPSS can be achieved by both reducing and increasing the size of DMU. %%%According to Theorem 3 in \textcolor{blue}{Podinovski (2004a)}, the GRS of DMU$_o$ in FDH models can be defined, based upon comparing $\theta_o^{FDH_C}, \theta_o^{FDH_{NI}},$ and $\theta_o^{FDH_{ND}}$ values, as follows.
%%%%------------------------------
\begin{definition}\label{KVE}\textcolor{blue}{(Podinovski 2004a)}: Let DMU$_o=(x_o,y_o)$ be an FDH$_V$-efficient unit. Then we have,\\
(i) G-CRS prevails at DMU$_o$ if and only if $\theta_o^{FDH_C}=\theta_o^{FDH_{NI}}=\theta_o^{FDH_{ND}}=1.$\\
(i) G-SCRS prevails at DMU$_o$ if and only if $\theta_o^{FDH_C}=\theta_o^{FDH_{NI}}=\theta_o^{FDH_{ND}}<1.$\\
(iii) G-IRS prevails at DMU$_o$ if and only if $\theta_o^{FDH_{ND}}>\theta_o^{FDH_C}=\theta_o^{FDH_{NI}}.$\\
(iv) G-DRS prevails at DMU$_o$ if and only if $\theta_o^{FDH_{NI}}>\theta_o^{FDH_C}=\theta_o^{FDH_{ND}}.$
\end{definition}
%=------------------------------
%%%%%%%%%%%%%%%%%%%%%%%%%%%%%%%%%%%%%%%%%%%%%%%%%%%%%%%%%%%%%%%%%%%%%%%%%%%

\section{Scale elasticity in FDH models: A drawback}
Let DMU$_o=(x_o,y_o)$ be the unit under consideration. The following RF has been used in the literature for defining Scale Elasticity (SE) in convex production technologies; see e.g. p. 150 in \textcolor{blue}{Podinovski et al. (2009)} and p. 1745 in \textcolor{blue}{Podinovski and F${\o}$rsund (2010)}:

\begin{equation}\label{function1}\beta_o(\alpha):=\max\{\beta~\vert~(\alpha x_o, \beta
y_o)\in T\},\end{equation}
in which $T$ is the corresponding production technology. In fact, $\beta_o(\alpha)$ gives the maximum proportion of the output
vector $y_o$, feasible in technology $T$ for the given
input vector $\alpha x_o$. The domain of this function is
%----------------------
\begin{equation}\label{lamlam}
\Lambda:=\{\alpha\in\mathbb{R}~\vert~\exists \beta\in
\mathbb{R}~\textmd{ s.t. }~(\alpha x_o, \beta y_o)\in
T\}. \end{equation}
%------------------------
Let $\varphi_o^T$ denote the output radial efficiency of unit $(x_o,y_o)$ under PPS $T$. The scale elasticity $\varepsilon(x,y)$ at any production point $(x,y)=(\alpha x_o,\beta_o(\alpha)y_o)$ is defined as the ratio of its marginal productivity $\beta_o^{'}(\alpha)$ (if it exists) to its average productivity $\frac{\beta_o(\alpha)}{\alpha}:$
\begin{equation}\label{SE}
\varepsilon(x,y)=\beta_o^{'}(\alpha)\frac{\alpha}{\beta_o(\alpha)}. \end{equation}
In particular at the unit $(x_o,y_o)$ with $\varphi_o^T=1$, we have $\alpha=1$, $\beta_o(\alpha)=\beta_o(1)=\varphi_o^T=1$, and hence
\begin{equation}\label{SE1}
\varepsilon(x_o,y_o)=\beta_o^{'}(1). \end{equation}
If $\beta_o^{'}(1)$ does not exists, then one-sided scale elasticities are defined as
\begin{equation}\label{SE2}
\varepsilon_+(x_o,y_o)=\beta_{o+}^{'}(1),\end{equation}
\begin{equation}\label{SE3}
\varepsilon_-(x_o,y_o)=\beta_{o-}^{'}(1).\end{equation}
\textcolor{blue}{Podinovski et al. (2009)} derived a simple formula to calculate $\beta_o^{'}(1),~\beta_{o+}^{'}(1)$, and $\beta_{o-}^{'}(1)$ under convexity of $T$. The following simple example shows that Eqs. (\ref{SE1})-(\ref{SE3}) cannot be used to define SE in nonconvex FDH technologies.

\begin{example}\label{EXMAMPLE}
Consider a set of four DMUs $A=(1,2),~B=(3,4),~C=(5,5),$ and $D=(6,13)$. Each DMU utilizes a single-input to produce a
single-output.\vspace{2cm}\\

\begin{picture}(400,300)(140,-20)

\put(180,10){\vector(1,0){310} } \put(190,0){\vector(0,1){300} }
\put(480,-2){Input}\put(180,310){Output}

%\put(178,0){O}\put(180,22){1}
\put(180,42){2}\put(180,82){4}\put(180,102){5}\put(173,262){13}

\put(228,0){1}
\put(308,0){3} \put(388,0){5}\put(428,0){6}

\put(230,45){\circle{2}}\put(310,85){\circle{2}}\put(390,105){\circle{2}}\put(430,265){\circle{2}}

\put(220,42){A}\put(298,83){B}\put(383,109){C}\put(420,268){D}

\put(230,10){\line(0,1){34}}
\put(230,45){\line(1,0){40}}
%\put(270,25){\line(0,1){21}}

\put(270,45){\line(1,0){40}} \put(310,45){\line(0,1){40}}
\put(310,85){\line(1,0){80}}
\put(390,85){\line(0,1){20}}
\put(390,105){\line(1,0){40}}
\put(430,105){\line(0,1){160}}
\put(430,265){\line(1,0){30}}
\put(462,265){$\ldots$}

\put(230,10){\circle{1}}\put(310,10){\circle{1}}\put(390,10){\circle{1}}
\put(430,10){\circle{1}}

\put(190,45){\circle{1}}\put(190,85){\circle{1}}\put(190,105){\circle{1}}\put(190,265){\circle{1}}

\end{picture}\\
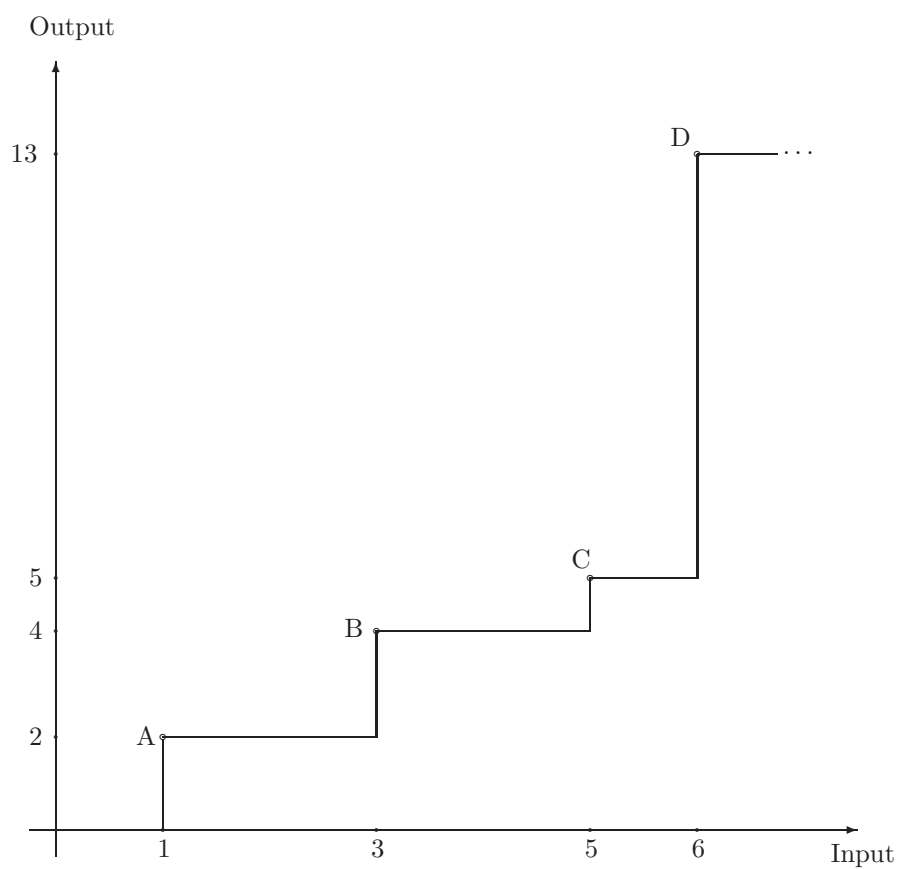
$$\textmd{ Figure 1. The FDH$_V$ production possibility set in Example \ref{EXMAMPLE}.}$$

The frontier of FDH$_V$ PPS can be seen
in Fig. 1. Consider DMU$_B$. From Fig. 1, it is seen that
$$\beta_B(\alpha)=\left\{\begin{array}{llll}
\frac{1}{2},&&&\frac{1}{3}<\alpha<1\vspace{1mm}\\
1,&&&1\leq \alpha<\frac{5}{3}
\end{array}\right.$$
Hence, $\beta_{B+}^{'}(1)=0$ and $\beta_{B-}^{'}(1)=\infty$.
\end{example}
%====================================================================

Indeed, as seen in the above example, $\beta_{o+}^{'}(1)$ and $\beta_{o-}^{'}(1)$ are always belong to $\{0,\infty\}$ if $T=FDH_V$ and DMU$_o$ is FDH$_V$-efficient. It is due to the piecewise constant structure of the frontier of FDH$_V$. Hence, Eqs. (\ref{SE1})-(\ref{SE3}) fail in defining SE in nonconvex FDH$_V$ technology.

\section{SE Counterpart}

GRS is an indicator of the direction in which the MPSS of an efficient DMU is achieved. When Increasing Returns to Scale (IRS) prevails at some DMU under a convex technology, proportional increase in its inputs leads to a greater increase in its outputs until getting MPSS. But in nonconvex technologies the matter is different. When G-IRS prevails at output-oriented efficient DMU$_o$, this unit has to increase its scale of operations to achieve MPSS (\textcolor{blue}{Podinovski 2004a}), but due to the non-convexity of the production possibility set, all increments in inputs may not lead to more increment in outputs (even under G-IRS).\\
For example, in Fig 1, G-IRS prevails at DMU$_A$ (due to Theorem 3 in  \textcolor{blue}{Podinovski 2004a}). If one increases $x_A$ from 1 to 3, then $y_A$ will increase from 2 to 4. The proportion of increase in input is equal to 3, while it is 2 for output\footnote{This property can be interpreted as non-monotonicity of RAP; see \textcolor{blue}{Podinovski (2004a,b)}}. It shows that the GRS status of units gives insufficient information to the manager(s). So, we need another criteria to decide about the amount of increment/decrement in the operation of some unit. In fact, $\lq\lq$the GRS classification suggest the types of resizing which the efficient unit should implement in order to achieve the global maximum of the average productivity"; see p. 239 in \textcolor{blue}{Podinovski (2004a)}. But GRS classification does not say anything about the quantity of the ratio among inputs and outputs in the path of getting MPSS. SE measure overcome this pitfall in convex DEA, but as seen in the previous section the traditional SE notion does not work under nonconvex FDH technologies. So, to have an insight about the ratio among proportion of outputs and proportion of inputs we need a counterpart of SE which works properly under nonconvex technologies. The rest of this section is devoted to preparing this material.

Let DMU$_o=(x_o,y_o)$ be the unit under consideration. The RF under FDH$_V$ PPS is defined as (\ref{function1}), and its domain is the set presented in (\ref{lamlam}).

It can be shown that $\Lambda=[\underline{\alpha}, +\infty)$ for some $\underline{\alpha}\leq 1$.

Here, we consider the right derivative, because G-IRS prevails at DMU$_o$. The remain cases can be analyzed analogously. Under FDH$_V$ PPS, the RF, $\beta_o(.)$, is a stepwise noncontinuous function with $\beta^{'}_{o+}(1)=0.$  Now, we are ready to define a SE counterpart to have a quantitative
measure of the strength of the GRS for DMUs under FDH technologies. As local SE, utilizing derivative of the response function at one, on the right of DMU$_o$ is zero, we define a maximum incremental ratio among inputs and outputs as follows:

\begin{equation}\label{rast+}
\sigma_o^+:=\displaystyle\max_{\alpha>1}\frac{\beta_o(\alpha)-1}{\alpha-1},
\end{equation}

The quantity $\sigma_o^+$ can be interpreted as a finite-difference value instead of the derivative $\beta_{o+}^{'}(1)=\displaystyle\lim_{\alpha\longrightarrow 1^+}\frac{\beta_o(\alpha)-1}{\alpha-1}$. In the rest of this section, we show that $\sigma_o^+$ can be obtained via some simple ratios. In Section 6, it is shown that the newly-defined quantity $\sigma_o^+$ is consistent with GRS.

Indeed, $\sigma_o^+$ is calculated by polynomial-time Procedure 1, invoking the following main theorem.
%======================================================================
\begin{theorem}\label{++++} Assume that DMU$_o=(x_o,y_o)$ is the unit under consideration, it is FDH$_V$-efficient and G-IRS prevails at this unit. Define
$$\alpha_{jo}:=\max_{i}\frac{x_{ij}}{x_{io}},~~~\beta_{jo}:=\displaystyle\min_{r}\frac{y_{rj}}{y_{ro}},~~~\Pi:=\{j\in J:~\alpha_{jo}>1\}.$$ Then:\\
(i) $\Pi\neq \emptyset$.\\
(ii) $\sigma_o^+=\displaystyle\max_{j\in \Pi}\frac{\beta_{jo}-1}{\alpha_{jo}-1}$.\\
\end{theorem}
%======================================================================

%=======================================================
%----------------------------------------------------
$$\begin{tabular}{|ll|}
  \hline
  \textbf{\small Procedure 1.}&\\\hline
\textit{Step 0.}& Let FDH$_V$-efficient $DMU_o=(x_o,y_o)$ with G-IRS status be given.\\
%------------------------------------------
\textit{Step 1.}& For $j=1$ to $n$ set
$\alpha_{jo}=\displaystyle\max_{i}\frac{x_{ij}}{x_{io}}$ and
$\beta_{jo}=\displaystyle\min_{r}\frac{y_{rj}}{y_{ro}}$. \\
%------------------------------------------
\textit{Step 2.}& Set $\Pi=\{j\in J:~\alpha_{jo}>1\}$. Then $\sigma_o^+=\displaystyle\max_{j\in \Pi}\frac{\beta_{jo}-1}{\alpha_{jo}-1}$.\\
 \hline
\end{tabular}$$\\
%----------------------------------------------------
%======================================================================

\begin{remark}
Although in our discussion in this section we focused on the G-IRS
case with $\alpha>1$ (increasing in inputs) and $\sigma^+_o$, the
G-DRS case with $\alpha<1$ (decreasing in inputs) can be discussed analogously. In G-DRS case, SE counterpart is defined as
\begin{equation}\label{rast-}
\sigma_o^-:=\displaystyle\min_{\underline{\alpha}\leq \alpha<1}\frac{\beta_o(\alpha)-1}{\alpha-1},
\end{equation}
and it is called the minimum decremental ratio among inputs and outputs. When G-SCRS prevails, the DMU$_o$ does not operate at MPSS but can choose whether to increase or reduce the scale of its operations in order to achieve MPSS. In such a case, two quantities $\sigma_o^+$ and $\sigma_o^-$ can be used based upon the decision of DM to increase or decrease the scale of operations, respectively.
\end{remark}

\section{One-sided RTS}

There are many papers in the DEA literature discussing the theory
and applications of RTS, see, e.g., \textcolor{blue}{Banker et
al. (2004)} for a review. One of the basic and useful PPS-based
definitions of RTS in DEA models has been introduced by
\textcolor{blue}{Banker et al. (1984),} and improved by
\textcolor{blue}{Tone (2001)} and
\textcolor{blue}{Soleimani-damaneh (2012)}. In this section, we
provide a one-sided counterpart of this definition in FDH
technologies. Establishing the relationships between one-sided RTS, GRS, and maximum/minimum incremental/decremental ratios is done in Section 6.

We start with explaining our motivation. Assume that G-IRS prevails at FDH$_V$-efficient DMU$_o$ under consideration. It means that this unit should expand its operation to achieve MPSS (maximum average productivity), though expanding the operation of a unit may not be possible in practice, due to budget limitations, physical restrictions, etc. Furthermore, increasing the scale of operations to some value before getting MPSS may not lead to a greater average productivity. In such a case, contracting the operation of the DMU might be more beneficial from average productivity standpoint. For example, in Fig. 1, consider DMU$_B$ with average productivity $\frac{4}{3}$. Assume that, due to the above-mentioned limitations, the maximum possible increase in input of this DMU is 2. Then DMU$_C$ is gotten whose average productivity equals one (less than that of DMU$_B$). On the other hand, decreasing the input of DMU$_B$ from 3 to 1 leads to DMU$_A$ with better average productivity. This fact guides us to examine the behavior of DMUs in both sides. To this end, we define one-sided RTS notion as follows.

Hereafter, for a given set $A$, the notations $intA$ and $A^c$ stand for the interior and the complement of $A$,
respectively.
%--------------------

\begin{definition} \label{DEF2}
Let $DMU_o=(x_o,y_o)$ be an FDH$_V$-efficient unit, and $$z_\delta=(\delta x_o, \delta y_o).$$ Then\\
(i) Right-IRS prevails at $DMU_o$ if $z_\delta\in int T^{FDH_V}$ for some $\delta>1.$\\
(ii) Right-DRS prevails at $DMU_o$ if $z_\delta\in (T^{FDH_V})^c$ for each $\delta>1.$\\
(iii) Right-CRS prevails at $DMU_o$ if neither Right-IRS nor Right-DRS prevails at this unit.
\end{definition}

For example, in Fig. 2, Right-IRS prevails at DMUs $A,~B,~D,$ and $F$, while Right-DRS prevails at DMUs $G$ and $E$; and Right-CRS prevails at DMU $C$.
%==================================

\begin{definition} \label{DEF22}
Let $DMU_o=(x_o,y_o)$ be an FDH$_V$-efficient unit. Then\\
(i) Left-IRS prevails at $DMU_o$ if $z_{\delta}\in (T^{FDH_V})^c$ for each $0<\delta<1.$\\
(ii) Left-DRS prevails at $DMU_o$ if $z_{\delta}\in int T^{FDH_V}$ for some $0<\delta<1$. \\
(iii) Left-CRS prevails at $DMU_o$ if neither Left-IRS nor Left-DRS prevails at this unit.
\end{definition}
%==================================

For example, in Fig. 2, Left-IRS prevails at DMUs $A$ and $C$, while Left-DRS prevails at DMUs $B,~D,~F,$ and $G$; and Left-CRS prevails at DMU $E$.
%==================================

The following lemma helps us in the sequel; its proof can be
found in \textcolor{blue}{Soleimani-damaneh and Mostafaee (2015).}

\begin{lemma}\label{LMA} (Lemma 4.1 in \textcolor{blue}{Soleimani-damaneh and Mostafaee, 2015}): Consider the following set:
$$T^0=\{(x,y) ~\vert~ \sum_{j\in J}\lambda_jx_j<x,~ \sum_{j\in J}\lambda_jy_j>y>0,~ \sum_{j\in J}\lambda_j=1,~ \lambda\in(\{0,1\})^n\}.$$ Then we have $int T^{FDH_V}=T^0.$
\end{lemma}

The following theorem provides a necessary and sufficient
condition for checking the right-RTS status of $DMU_o=(x_o,
y_o)$. The variables of systems (\ref{EQNA})-(\ref{EQND}) in Theorems \ref{THMB} and \ref{THMC} are $\lambda_1,\lambda_2,\ldots,\lambda_n,\delta.$
Part (i) of the theorem results from Lemma \ref{LMA} and Definition \ref{DEF2}. Parts (ii) and (iii) are straightforward due to Definition \ref{DEF2}.

%=======================================================
\begin{theorem}\label{THMB}
Let $DMU_o=(x_o,y_o)$ be an FDH$_V$-efficient unit. Then\\
(i) Right-IRS prevails at $DMU_o$ if and only if the following system has some solution:
\begin{equation}\label{EQNA}
\left\{\begin{array}{ll}
\displaystyle\sum_{j\in J}\lambda_jx_j<\delta x_o,&\\
\displaystyle\sum_{j\in J}\lambda_jy_j>\delta y_o,&\\
\displaystyle \sum_{j\in J}\lambda_j=1,&\\
\lambda_j\in \{0,1\},~\delta>1, &j=1,\ldots,n.
\end{array}\right.
\end{equation}
(ii) Right-DRS prevails at $DMU_o$ if and only if the following system has no solution:\\
\begin{equation}\label{EQNB}
\left\{\begin{array}{ll}
\displaystyle\sum_{j\in J}\lambda_jx_j\leq\delta x_o,&\\
\displaystyle\sum_{j\in J}\lambda_jy_j\geq \delta y_o,&\\
\displaystyle \sum_{j\in J}\lambda_j=1,&\\
\lambda_j\in \{0,1\},~\delta>1, &j=1,\ldots,n.\\
\end{array}\right.
\end{equation}
(iii) Right-CRS prevails at $DMU_o$ if and only if System (\ref{EQNA}) does not have any solution and System (\ref{EQNB}) has some solution.
\end{theorem}
%=======================================================

Theorem \ref{THMC} results from Lemma \ref{LMA} and Definition \ref{DEF22}.

\begin{theorem}\label{THMC}
Let $DMU_o=(x_o,y_o)$ be an FDH$_V$-efficient unit. Then\\
(i) Left-IRS prevails at $DMU_o$ if and only if the following system has no solution:\\
\begin{equation}\label{EQNC}
\left\{\begin{array}{ll}
\displaystyle\sum_{j\in J}\lambda_jx_j\leq \delta x_o,&\\
\displaystyle\sum_{j\in J}\lambda_jy_j\geq \delta y_o,&\\
\displaystyle \sum_{j\in J}\lambda_j=1,&\\
\lambda_j\in \{0,1\},~\delta\in (0,1), &j=1,\ldots,n.\\
\end{array}\right.
\end{equation}
(ii) Left-DRS prevails at $DMU_o$ if and only if the following system has some solution:\\
\begin{equation}\label{EQND}
\left\{\begin{array}{ll}
\displaystyle\sum_{j\in J}\lambda_jx_j<\delta x_o,&\\
\displaystyle\sum_{j\in J}\lambda_jy_j>\delta y_o,&\\
\displaystyle \sum_{j\in J}\lambda_j=1,&\\
\lambda_j\in \{0,1\},~\delta\in (0,1), &j=1,\ldots,n.
\end{array}\right.
\end{equation}
(iii) Left-CRS prevails at $DMU_o$ if and only if System (\ref{EQNC}) has some solution and System (\ref{EQND}) has no solution.
\end{theorem}
%=====================================================

Theorems \ref{THMQ} and \ref{THM5} provide some ratio-based polynomial-time tests to obtain the one-sided RTS class of the FDH$_V$-efficient unit under consideration DMU$_o$. These theorems result from Theorems \ref{THMB} and \ref{THMC}.
The proof of Theorem \ref{THMQ} is given in the appendix. The proof of Theorem \ref{THM5} is similar to that of Theorem
\ref{THMQ}, and is hence omitted.

%=======================================
\begin{theorem} \label{THMQ}
Let $DMU_o=(x_o,y_o)$ be an FDH$_V$-efficient unit. Define
$$\alpha_{jo}:=\max_{i}\frac{x_{ij}}{x_{io}},~~~\beta_{jo}:=\displaystyle\min_{r}\frac{y_{rj}}{y_{ro}}.$$ Then we have\\
(i) Right-IRS prevails at DMU$_o$ if and only if there exists some $j\in J$ such that $\beta_{jo}>1$ and $\alpha_{jo}<\beta_{jo}.$\\
(ii) Right-DRS prevails at DMU$_o$ if and only if $\beta_{jo}\leq1$ or $\beta_{jo}<\alpha_{jo}$ for each $j\in J.$\\
(iii) Right-CRS prevails at DMU$_o$ if and only if both the following conditions hold:\\
$~~~~~~~~~~$(iii-a): $\alpha_{jo}\geq \beta_{jo}$ or $\beta_{jo}\leq 1$ for each $j\in J$\\
$~~~~~~~~~~~$(iii-b): $\beta_{jo}> 1$ and $\beta_{jo}\geq \alpha_{jo}$ for some $j\in J.$
\end{theorem}
%-------------

%==========================================================
\begin{theorem}\label{THM5}
Let $DMU_o=(x_o,y_o)$ be an FDH$_V$-efficient unit. Let $\alpha_{jo}$ and $\beta_{jo}$ be as defined in Theorem \ref{THMQ}. Then we have\\
(i) Left-IRS prevails at $DMU_o$ if and only if $\alpha_{jo}\geq1$ or $\beta_{jo}<\alpha_{jo}$ for each $j\in J.$\\
(ii) Left-DRS prevails at $DMU_o$ if and only if there exists some $j\in J$ such that $\alpha_{jo}<1$ and $\alpha_{jo}<\beta_{jo}.$\\
(iii) Left-CRS prevails at $DMU_o$ if and only if both the following conditions hold:\\
$~~~~~~~~~$(iii-a): $\alpha_{jo}<1$ and $\beta_{jo}\geq \alpha_{jo}$ for some $j\in J$.\\
$~~~~~~~~~~$(iii-b): $\alpha_{jo}\geq 1$ or $\alpha_{jo}\geq \beta_{jo}$ for each $j\in J$.
\end{theorem}
%=======================================================

Although, invoking Theorems \ref{THMQ} and \ref{THM5}, a polynomial-time procedure can be provided to determine the one-sided RTS class of the unit under consideration DMU$_o$, the following theorem shows that it can be derived from Procedure 1 as well. We close this section with this important result which presents the relationship between one-sided RTS and maximum/minimum incremental/decremental ratios.
%=============================
\begin{theorem}\label{ipm0}
Let $DMU_o=(x_o,y_o)$ be an FDH$_V$-efficient unit. Then,\\
(i) Right-IRS prevails at DMU$_o$ if and only if $\sigma_o^+>1.$\\
(ii) Right-DRS prevails at DMU$_o$ if and only if $\sigma_o^+<1.$\\
(iii) Right-CRS prevails at DMU$_o$ if and only if $\sigma_o^+=1.$\\
(iv) Left-IRS prevails at DMU$_o$ if and only if $\sigma_o^->1.$\\
(v) Left-DRS prevails at DMU$_o$ if and only if $\sigma_o^-<1.$\\
(vi) Left-CRS prevails at DMU$_o$ if and only if $\sigma_o^-=1.$\\
\end{theorem}
%========================================================

\section{One-sided RTS, GRS, and maximum/minimum ratios}

Theorem \ref{ipm1} addresses a connection between one-sided RTS, GRS, and maximum/minimum ratios. Here, right-NIRS stands for prevailing right-CRS or right-DRS. Also, left-NDRS means prevailing left-CRS or left-IRS.
\begin{theorem}\label{ipm1}
Let $DMU_o=(x_o,y_o)$ be an FDH$_V$-efficient unit. Then,\\
(a) G-IRS $\Longrightarrow$ Right-IRS $\Longleftrightarrow \sigma_o^+>1.$\\
(b) G-DRS $\Longrightarrow$ Left-DRS $\Longleftrightarrow \sigma_o^-<1.$\\
(c) G-CRS $\Longrightarrow$ Right-NIRS and Left-NDRS  $\Longleftrightarrow \sigma_o^+\leq 1$ and $\sigma_o^-\geq 1$.\\
(d) G-SCRS $\Longrightarrow$ Right-IRS and Left-DRS  $\Longleftrightarrow \sigma_o^+>1$ and $\sigma_o^-<1$.\\
\end{theorem}

Part (d) of the above theorem reveals a difference between G-CRS and G-SCRS statuses. This shows that under G-SCRS none of the one-sided RTS positions is constant.

%%%Notice that in contrast to convex technologies (\textcolor{blue}{Jahanshahloo et al. 2005}), in FDH technologies Right-IRS (resp. Left-DRS) does not imply G-NDRS (resp. G-NIRS) necessarily. For instance, in Example \ref{EXMA}, DMU $F$ has right-IRS status and DMU $B$ has left-DRS status. At these two units G-DRS and G-IRS prevails, respectively.

\section{Conclusions}

In this paper, two counterparts for one-sided SE in FDH models, as a popular class of nonconvex DEA models, have been defined and investigated. Using a response function, two measures, called maximum incremental and minimum decremental ratios among inputs and outputs, have been defined. A theorem has been established to calculate the ratios. In the second part of the paper, one-sided RTS has been defined followed by some theoretical results to determine it. Final theorem of the paper proves the relationship between maximum incremental and minimum decremental ratios, one-sided RTS, and GRS. Resulting from theoretical results, a polynomial-time procedure has been sketched which is able to calculate maximum incremental and minimum decremental ratios, and to determine the one-sided RTS of DMUs.

\section{Appendix: Proofs of the main results}
\textbf{Theorem \ref{++++}}. Assume that DMU$_o=(x_o,y_o)$ is the unit under consideration, it is FDH$_V$-efficient and G-IRS prevails at this unit. Define
$$\alpha_{jo}:=\max_{i}\frac{x_{ij}}{x_{io}},~~~\beta_{jo}:=\displaystyle\min_{r}\frac{y_{rj}}{y_{ro}},~~~\Pi:=\{j\in J:~\alpha_{jo}>1\}.$$ Then:\\
(i) $\Pi\neq \emptyset$.\\
(ii) $\sigma_o^+=\displaystyle\max_{j\in \Pi}\frac{\beta_{jo}-1}{\alpha_{jo}-1}$.\\
%----------------------
\textbf{Proof}. (i): Define \begin{equation}\label{mineral1}
S^1:=\{j\in J: \displaystyle \beta_{jo}>1\},~~~~~~~~~~~~~
\end{equation}
Since G-IRS prevails at DMU$_o$, we have $\theta_o^{FDH_C}<1$, leading to
$\varphi_o^{FDH_{NI}}=\varphi_o^{FDH_C}>1$. Hence, there
exists some $j_0\in J$ and some $\lambda_{j_0}\in (0,1]$ satisfying
$$\frac{y_{rj_0}}{y_{ro}}\geq \frac{\lambda_{j_0}y_{rj_0}}{y_{ro}}\geq \varphi_o^{FDH_{NI}}>1,~~\forall r.$$
This implies $j_0\in S^1$.

So, $S^1\neq \emptyset$, and to prove the nonemptyness of $\Pi$,
it is sufficient to show that $S^1\subseteq \Pi$. Let $j\in S^1$.
Then $y_j>y_o$. If $j\notin \Pi$, then $\alpha_{jo}\leq 1$ which
leads to $x_j\leq x_o$. Hence, we get $x_j\leq x_o$, $y_j\geq
y_o$, and $(x_j,y_j)\neq (x_o,y_o)$. This contradicts the
FDH$_V$-efficiency of DMU$_o$; and the proof of this part is completed. \\
(ii): It is not difficult to see that
\begin{equation}\label{Eade}
\displaystyle\begin{array}{lll}
\sigma_o^+=&\displaystyle\max_{\alpha,\beta}& \frac{\beta-1}{\alpha-1}\\
&s.t.&\displaystyle\sum_{j\in J} \lambda_jx_j\leq \alpha x_o,\\
&&\displaystyle\sum_{j\in J} \lambda_j y_j \geq \beta y_o,\\
&&\displaystyle\sum_{j\in J}\lambda_j=1,\\
&&\lambda_j \in \{0,1\};~j\in J,\\
&&\alpha>1,~\beta\geq 0.\\
\end{array}
\end{equation}
Let $\Upsilon^*:=(\lambda_k^*=1,~\lambda_j^*=0;~~j\in
J\backslash\{k\},~\beta^*,~\alpha^*)$ be an optimal solution to
Problem (\ref{Eade}). Then $\frac{x_{ik}}{x_{io}}\leq \alpha^*$
for each $i$, and $\frac{y_{rk}}{y_{ro}}\geq\beta^*$ for each
$r$. Because of the optimality of $\Upsilon^*$ and due to the
objective function of Problem (\ref{Eade}), we get
%------
$$1<\alpha^*=\displaystyle\max_i\frac{x_{ik}}{x_{io}}=\alpha_{ko},~~\textmd{ and }~~\beta^*=\displaystyle\min_r\frac{y_{rk}}{y_{ro}}=\beta_{ko}.$$
%-----------
These imply that $k\in \Pi$ and $k\neq o$. Furthermore,
$$\sigma_o^+=\frac{\beta_{ko}-1}{\alpha_{ko}-1}.$$
To complete the proof, we should show that
$$\frac{\beta_{ko}-1}{\alpha_{ko}-1}\geq
\frac{\beta_{jo}-1}{\alpha_{jo}-1},~~\forall j\in \Pi.$$ By
contradiction, assume that
$\frac{\beta_{ko}-1}{\alpha_{ko}-1}<\frac{\beta_{po}-1}{\alpha_{po}-1}$
for some $p\in \Pi.$ Then $(\lambda_p=1,~\lambda_j=0;~~j\in
J\backslash\{p\},~\beta=\beta_{po},~\alpha=\alpha_{po})$ is a
feasible solution to Problem (\ref{Eade}) with greater objective
function value than $\sigma_o^+$. This contradicts the optimality
of $\Upsilon^*$ and the proof is completed.\qed$~$\\
%=======================================
\textbf{Theorem \ref{THMQ}}.
Let $DMU_o=(x_o,y_o)$ be an FDH$_V$-efficient unit. Define
$$\alpha_{jo}:=\max_{i}\frac{x_{ij}}{x_{io}},~~~\beta_{jo}:=\displaystyle\min_{r}\frac{y_{rj}}{y_{ro}}.$$ Then we have\\
(i) Right-IRS prevails at DMU$_o$ if and only if there exists some $j\in J$ such that $\beta_{jo}>1$ and $\alpha_{jo}<\beta_{jo}.$\\
(ii) Right-DRS prevails at DMU$_o$ if and only if $\beta_{jo}\leq1$ or $\beta_{jo}<\alpha_{jo}$ for each $j\in J.$\\
(iii) Right-CRS prevails at DMU$_o$ if and only if both the following conditions hold:\\
$~~~~~~~~~~$(iii-a): $\alpha_{jo}\geq \beta_{jo}$ or $\beta_{jo}\leq 1$ for each $j\in J$\\
$~~~~~~~~~~$(iii-b): $\beta_{jo}> 1$ and $\beta_{jo}\geq \alpha_{jo}$ for some $j\in J.$\\
%-------------
\textbf{Proof}. (i): By Theorem \ref{THMB}, Right-IRS prevails at DMU$_o$ if and only if there exists some $j\in J$ such that the following system has a solution $\delta$:
\begin{equation}
\left\{\begin{array}{ll}
\displaystyle x_{ij}<\delta x_{io},&i=1,\ldots,m,\\
\displaystyle y_{rj}>\delta y_{ro},&r=1,\ldots,s,\\
~\delta>1.

\end{array}\right.
\end{equation}
It holds if and only if there exists some
$j \in J$  and some $\delta>1$ satisfying
\begin{equation}
\alpha_{jo}=\displaystyle\max_{i}\{\frac{x_{ij}}{x_{io}}\}<\delta<\min_{r}\{\frac{y_{rj}}{y_{ro}}\}=\beta_{jo}
\end{equation}
It holds if and only if there exists some $j \in J$ such that $\beta_{jo}>1$ and $\alpha_{jo}<\beta_{jo}$. This completes the proof of Part (i).\\
%-------------
(ii): By Theorem \ref{THMB}, Right-DRS does not prevails at DMU$_o$ if and only if there exists some $j\in J$ and some $\delta>0$ satisfying
\begin{equation}\label{EQUATIONE}
\left\{\begin{array}{l}
\displaystyle x_j\leq\delta x_o,\\
\displaystyle y_j\geq \delta y_o,\\
~\delta>1. \\
\end{array}\right.
\end{equation}
Therefore,
Right-DRS does not prevails at DMU$_o$ if and only if there exists some $j\in J$ and some $\delta>1$ satisfying
$\alpha_{jo}\leq \delta \leq \beta_{jo}.$
It holds if and only if there exists some $j\in J$ such that $\beta_{jo}>1$ and $\alpha_{jo}\leq \beta_{jo}.$ Hence, Right-DRS prevails at $DMU_o$ if and only if $\beta_{jo}\leq 1$ or $\alpha_{jo}> \beta_{jo}$ for each $j \in J.$ This completes the proof of Part (ii).\\
(iii): Part (iii) results from parts (i) and (ii).\qed$~$\\
%-----------------------------------
\textbf{Theorem \ref{ipm0}}.
Let $DMU_o=(x_o,y_o)$ be an FDH$_V$-efficient unit. Then,\\
(i) Right-IRS prevails at DMU$_o$ if and only if $\sigma_o^+>1.$\\
(ii) Right-DRS prevails at DMU$_o$ if and only if $\sigma_o^+<1.$\\
(iii) Right-CRS prevails at DMU$_o$ if and only if $\sigma_o^+=1.$\\
(iv) Left-IRS prevails at DMU$_o$ if and only if $\sigma_o^->1.$\\
(v) Left-DRS prevails at DMU$_o$ if and only if $\sigma_o^-<1.$\\
(vi) Left-CRS prevails at DMU$_o$ if and only if $\sigma_o^-=1.$\\
%----
\textbf{Proof}. (i): If Right-IRS prevails at DMU$_o$, then by Definition \ref{DEF2}, there exists some $\delta>1$ with $z_\delta\in int T^{FDH_V}.$ Therefore, there exists some scalar $\varepsilon\in (0,\delta)$ such that $\delta-\varepsilon>1$ and $$z_\delta+\varepsilon(-x_o,y_0)=((\delta-\varepsilon)x_o,(\delta+\varepsilon)y_o)\in T^{FDH_V}.$$ This implies $\beta_o(\delta-\varepsilon)\geq \delta+\varepsilon.$ Hence,
 %---------
$$\sigma_o^+=\displaystyle\max_{\alpha>1}
\frac{\beta_o(\alpha)-1}{\alpha-1}\geq
\frac{\beta_o(\delta-\varepsilon)-1}{\delta-\varepsilon-1}\geq
\frac{\delta+\varepsilon-1}{\delta-\varepsilon-1}>1.$$
%-----------------------------

Conversely, assume that $\sigma_o^+>1$. Then, there exists some $\alpha>1$ such that $\frac{\beta_o(\alpha)-1}{\alpha-1}>1,$ and $(\alpha x_o, \beta_o(\alpha)y_o)\in T^{FDH_V}.$ Set $$\delta=\frac{\beta_o(\alpha)+\alpha}{2}.$$ Then $\beta_o(\alpha)>\delta>\alpha>1$. Since $(\alpha x_o, \beta_o(\alpha)y_o)\in T^{FDH_V},$ there exists some $\lambda\in (\{0,1\})^n$ such that $\displaystyle\sum_{j\in  J}\lambda_j=1$ and
\begin{equation}
\displaystyle\sum_{j\in J}\lambda_jx_j\leq \alpha x_o<\delta x_o,
\end{equation}
\begin{equation}
\displaystyle\sum_{j\in J}\lambda_jy_j\geq \beta_o(\alpha) y_o>\delta y_o.
\end{equation}
These imply $z_\delta\in int T^{FDH_V},$ because of Lemma \ref{LMA}. Hence, Right-IRS prevails at DMU$_o$ according to Definition \ref{DEF2}. \\
%---------------------------------------------------------------
(ii): Assume that Right-DRS prevails at DMU$_o$. Then, by Definition \ref{DEF2}, $z_\delta\in (T^{FDH_V})^c$ for each $\delta>1.$ Therefore, according to the possibility (free disposability) axiom, we have $\beta_o(\delta)<\delta$ for each $\delta>1.$ Hence,
%-----
$$\sigma_o^+=\displaystyle\max_{\delta>1}\frac{\beta_o(\delta)-1}{\delta-1}<1.$$
%-------------------

Conversely, assume that $\sigma_o^+<1.$ Thus $\beta_o(\delta)<\delta$ for each $\delta>1.$ This means that $z_\delta \notin T^{FDH_V}$ for each $\delta>1,$ and the proof of part (ii) is completed.  \\
%6666666666666666666666666666
(iii)-(vi): Part (iii) results from parts (i) and (ii). The proof of parts (iv)-(vi) are similar to that parts (i)-(iii).\qed$~$\\
%====================
\textbf{Theorem \ref{ipm1}}.
Let $DMU_o=(x_o,y_o)$ be an FDH$_V$-efficient unit. Then,\\
(a) G-IRS $\Longrightarrow$ Right-IRS $\Longleftrightarrow \sigma_o^+>1.$\\
(b) G-DRS $\Longrightarrow$ Left-DRS $\Longleftrightarrow \sigma_o^-<1.$\\
(c) G-CRS $\Longrightarrow$ Right-NIRS and Left-NDRS  $\Longleftrightarrow \sigma_o^+\leq 1$ and $\sigma_o^-\geq 1$.\\
(d) G-SCRS $\Longrightarrow$ Right-IRS and Left-DRS  $\Longleftrightarrow \sigma_o^+>1$ and $\sigma_o^-<1$.\\
\textbf{Proof}. (a): Let G-IRS prevails at DMU$_o$. We have
$\varphi_o^{FDH_{NI}}>1$. Hence, there
exists some $j\in J$ and some $\lambda_j\in (0,1)$ satisfying
$$\lambda_jx_j\leq x_o,~~~\lambda_jy_j\geq \varphi_o^{FDH_{NI}}y_o.$$
Thus, there exists some $\varepsilon>0$ such that $\varphi_o^{FDH_{NI}}-\varepsilon>1$ and
$$\lambda_jx_j\leq x_o<(\varphi_o^{FDH_{NI}}-\varepsilon)x_o,~~~\lambda_jy_j\geq \varphi_o^{FDH_{NI}}y_o>(\varphi_o^{FDH_{NI}}-\varepsilon)y_o.$$
By setting $\delta:=\frac{\varphi_o^{FDH_{NI}}-\varepsilon}{\lambda_j}$, we have $\delta>1$ and $z_\delta\in intT^{FDH_V}$ according to Lemma \ref{LMA}. Hence, right-IRS prevails at DMU$_o$. Now the proof of part (a) is completed due to Theorem \ref{ipm0}.\\
(b)-(d): The proofs of these parts are similar to that Part (a) and are hence omitted.\qed

\end{document}